# Quantized Method Solution for Various Fluid-Solid Reaction Models

H. Delavari Amrei[*], E. Jamshidi, H. Ale Ebrahim

*Department of Chemical Engineering, Amirkabir University of Technology (Tehran Polytechnic), Petrochemical Center of Excellency, Tehran 15875-4413, Iran.*

[*]*Corresponding author. Tel.: +98 21 64543177; fax: +98 21 66405847.*

*E-mail address:Delavarih@aut.ac.ir (H. Delavari Amrei)*

**Abstract**

Fluid–solid reactions exist in many chemical and metallurgical process industries. Several models describe these reactions such as volume reaction model, grain model, random pore model and nucleation model. These models give two nonlinear coupled partial differential equations (CPDE) that must be solved numerically. A new approximate solution technique (quantized method) has been introduced for some of these models in recent years. In this work, the various fluid-solid reaction models with their quantized and numerical solutions have been discussed.

**Contents**

## 1. *Introduction*

Non-catalytic fluid–solid reactions are very important in many chemical and metallurgical processes such as metal oxides reduction [1-5] roasting of metallic sulfides [6, 7], adsorption of acid gases [8, 9], coal gasification [10, 11], activated carbon production [12-14] and catalyst regeneration [15, 16]. In order to interpret laboratory data on these systems and in design and scale up, the mathematical modeling of the single pellet reaction is very important. These systems have received considerable attention and lots of models and techniques for their solution are available in the specialized literature.

In fluid-solid reactions, the solid may be initially nonporous or porous. The reaction of nonporous solids can be divided into three types of geometrical groups [17]. For mathematical modeling of nonporous solid, the sharp interface model is used. This model is one of the earliest models used and its analytical solution is well described in standard textbooks on chemical reaction engineering [18-30]. Also, some commonly used models for a porous solid pellet consist of single pore model [31-37], volume reaction model [4, 30, 38-45], the grain model [46-62], random pore model [63-71] and nucleation model [72-75].

Modeling of Non-catalytic fluid–solid reactions for porous pellets leads to a pair of coupled partial differential equations and due to their complexity, these equations must be solved by numerical analysis or approximate methods. The common numerical method which have been in use are finite volume [47, 62], finite element [76, 77], finite difference such as Crank–Nicholson [40, 78], line method [79, 80], and collocation methods [81, 82].

The approximate methods have been extensively developed by a number of investigators and the literature is full of techniques for simplifying the associated mathematical and computational difficulties. The use of the cumulative concentration concept has been proposed independently by Dudukovic and Lamba [40] and Del Borghi et al. [83]. This transformation reduces the system of equation to a single nonlinear diffusion–reaction type equation in terms of a new variable namely the cumulative gas concentration. Sohn has proposed the law of additive reaction times for isothermal fluid-solid reaction, which produces relatively simple approximate solution to the governing CPDE. The law of additive reaction time, formulated and tested for shrinking-core model, grain models and random-pore models [58, 59, 74, 75, 80, 84- 87] but it has 50% error for some cases. A generalized Thiele Module approach has been used by Rmachandran [88]. Based on this method, an approximate solution obtained for some of the commonly encountered rate forms including the grain model. Also, Brem and Brouwers [89] presented a new description for the case of reaction rate of general order with respect to gas concentration. Despite some limitations, their formulation represents a pioneer attempt of finding an approximate solution which explicitly includes the combined effects of both non-linear kinetics and intrinsic solid surface development. Moreover, Marcos et al. [90] proposed a new approximate solution but their method involves some limitations.

Jamshidi and Ale Ebrahim proposed a new approach which is called the quantized method (QM). This method greatly permits to reduce the mathematical difficulties normally present in fluid–solid reaction problems. They illustrated the QM potential by applying it to some of gas–solid reaction models, including the grain model, half-order volume reaction model, nucleation model, and modified grain model [91-94]. Furthermore, Rafsanjani et al. [95] applied this method to propose a new mathematical solution for predicting char activation reactions. More recently, Shiravani et al. [96] illustrated the QM potential by applying it to several fluid–solid reactions and the ability of this method in the unsteady state fluid conservation problems. In another work, Gómez- Barea and Ollero [79] extended the range of application of the QM method. This was done by permitting to accommodate any general kinetics, i.e., $n$th-order, Michaelis–Menten, etc. and any explicitly given intrinsic behavior of the solid structure variation with reaction.

The above mentioned QM has shown that accuracy of results in comparison with other numerical methods for each model is good. Therefore, in this work, partial differential equations of the various fluid-solid reaction models with their QM solution have been reviewed. Also, this work presents the application of QM for the quasi-steady state random pore model, simultaneous (two gas with a solid) reactions, and gas-solid reactions in the packed beds with a good accuracy.

## *2. Mathematical models of single pellet reaction and their QM solutions*

In fluid-solid reactions with the porous solid, diffusion and chemical reaction occur simultaneously inside the pellet. Therefore, the reaction occurs in a diffuse zone rather than at a sharp boundary as for nonporous solids.

Consider the following fluid-solid reaction:

$$A(f) + \nu_B B(s) \rightarrow C(f) + D(s) \tag{1}$$

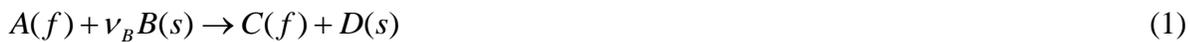

The following assumptions are applied generally for modeling of reactions in porous pellets:

(1) The reaction is irreversible and first-order with respect to the fluid reactant. The most of power law or lagmuir-Hinshelwood kinetics could be approximated by first-order reactions [17].

(2) The bulk fluid concentration is constant at the single pellet systems.

(3) The overall pellet size is constant. In the porous pellets, the structural change causes porosity to change and the pellet size usually remains constant.

(4) The bulk flow effect is negligible and the system is assumed as equimolar counter-diffusion.

(5) The system is isothermal for relatively low heat of reactions.

### 2.1. Volume reaction model

When the solid is porous, the fluid can penetrate into the solid and the reaction may now be assumed to take place all over the volume of the pellet. The models which represent this situation are known as volume reaction models. The rate of reaction at the interior points would be of course generally lower than that at the surface due to diffusional gradients [21].

These reactions are first order with respect to the fluid and can be zero, half or first order with respect to the solid reactant. The dimensionless governing equations of this model in general form are as follows [97]:

$$\psi\phi_V^2 \frac{\partial a}{\partial \theta_V} = \frac{\partial^2 a}{\partial y^2} + \frac{F_P - 1}{y} \frac{\partial a}{\partial y} - \phi_V^2 ab^n \quad (2)$$

$$\frac{\partial b}{\partial \theta_V} = -ab^n \quad (3)$$

Boundary conditions for eq. (2) are Dirichlet form as:

$$y = 0; \quad \frac{\partial a}{\partial y} = 0 \quad (4)$$

$$y = 1; \quad a = 1 \quad (5)$$

And initial condition for eq. (3) is:

$$\theta_V = 0; \quad b = 1 \quad (6)$$

Also, conversion of solid at each time can be calculated by the following equation:

$$X(\theta) = 1 - F_P \int_0^1 y^{F_P - 1} b \, dy \quad (7)$$

The numerical solution for the above coupled partial differential equations lead to a large set of coupled linear equations. By considering "i" as position index and "j" as time index, in equations (2,3), a(j,i), a(j,i+1), a(j, i-1), b(j,i), b(j-1,i), ... are related to each other. Therefore much computational effort is needed to obtain the final results.

In the quantized method, in contrast to the numerical solution, it is assumed that for small increments of time a(j,i) and b(j,i) are independent of a(j,i±1), b(j±1,i), y(i±1) and θ(j±1). In other words, the variables "a", "b", "θ" and "y" are related only by their (j,i) states. They are independent in term of their (j-1,i) or (j,i-1) [91-94].

By using b(j-1,i) instead of b(j,i) in Eq.(2) as an approximation, on the right-hand side of this equation, a(i,j) is the only variable and the other term can be treated as a constant. Therefore, a modified Thiele modulus is defined which converts Eqs. (2), (3) independent for small increments of time. For example, this modified Thiele modulus for unsteady state first order volume reaction model is defined as follows [96]:

$$M(j,i) = \phi_V b(j-1,i)^{1/2} \tag{8}$$

Now it is possible to solve equations (2,3) analytically for small time and position increments. Considering rate of reaction, pellet shape, boundary conditions, and accumulation term quantity, there are different cases for volume reaction models. The accumulation term in the fluid conservation for steady state cases is several orders of magnitude smaller than the diffusion and reaction terms. Thus quasi-steady state assumption is used for gas-solid reactions. The reason is that in gas–solid reactions usually $C_A/C_B \leq 10^{-3}$ holds and the accumulation term is negligible ($\psi=0$). However, in liquid-solid reactions such as phosphoric acid production or leaching processes the accumulation term must be considered.

Some approximate [40, 74, 79, 83, 88] and numerical solutions [40, 76- 78,80, 82] exist in the literature for volume reaction models. The analytical solution for small increments of time in the quantized technique for different cases of volume reaction model is presented in Tables 1-3. As presented in these tables, for half order reaction there are two stages. In the first stage, chemical reaction and diffusion occur, during which the outer solid in the pellet reacts until there is no longer any unreacted solid at the outer boundary of the pellet. Then the second stage is started, during which the boundary between the outer layer and the unreacted inner zone moves towards the centre of the pellet (moving boundary problem) [42].

Table 1. The analytical solution for small increments of time in the quantized technique for first order volume reaction model [92, 96].

**First order volume reaction model (n=1)**

$M = \phi_V b$

**quasi-steady state**

*Slab pellet ($F_P=1$)*

$$a = \frac{\cosh(My)}{\cosh(M)}$$

$$b = \exp\left[-\frac{\cosh(My)}{\cosh(M)}\theta_V\right]$$

*Spherical pellet ($F_P=3$)*

$$a = \frac{\sinh(My)}{y\sinh(M)}$$

$$b = \exp\left[-\frac{\sinh(My)}{y\sinh(M)}\theta_V\right]$$

**Unsteady state**

*Slab pellet ($F_P=1$)*

$$a = \frac{\cosh(My)}{\cosh(M)} + \sum_{k=1}^{\infty}\frac{k\pi(-1)^k \cos(k\pi y/2)}{M^2+(k^2\pi^2/4)}\exp\left[-\frac{(M^2+(k^2\pi^2/4))\theta_V}{\psi\phi_V^2}\right]$$

$$b = \exp\left[2-\frac{\cosh(My)}{\cosh(M)}\theta_V - \sum_{k=1}^{\infty}\frac{\frac{k\pi(-1)^k \cos(k\pi y/2)}{M^2+(k^2\pi^2/4)}\frac{\psi\phi_V^2}{M^2+(k^2\pi^2/4)}}{\times\left(1-\exp\left[-\frac{(M^2+(k^2\pi^2/4))\theta_V}{\psi\phi_V^2}\right]\right)}\right]$$

*Spherical pellet ($F_P=3$)*

$$a = \frac{\sinh(My)}{y\sinh(M)} + \sum_{k=1}^{\infty}\frac{2k\pi(-1)^k \sin(k\pi y)}{(M^2+k^2\pi^2)y}\exp\left[-\frac{(M^2+k^2\pi^2)\theta_V}{\psi\phi_V^2}\right]$$

$$b = \exp\left[2-\frac{\sinh(My)}{y\sinh(M)}\theta_V - \sum_{k=1}^{\infty}\frac{\frac{2k\pi(-1)^k \sin(k\pi y)}{(M^2+k^2\pi^2)y}\frac{\psi\phi_V^2}{M^2+k^2\pi^2}}{\times\left(1-\exp\left[-\frac{(M^2+k^2\pi^2)\theta_V}{\psi\phi_V^2}\right]\right)}\right]$$

Table 2. The analytical solution for small increments of time in the quantized technique for half order volume reaction model (quasi-steady state) [92].

**Half order volume reaction model (n=1/2)**

$M = \phi_V b^{1/2}$

   **quasi-steady state**

     *First stage*

        <u>Slab pellet ($F_P=1$)</u>

$$a = \frac{\cosh(My)}{\cosh(M)}$$

$$b = \left(1 - a\frac{\theta_V}{2}\right)^2$$

        <u>Spherical pellet ($F_P=3$)</u>

$$a = \frac{\sinh(My)}{y\sinh(M)}$$

$$b = \left(1 - a\frac{\theta_V}{2}\right)^2$$

     Second stage

        <u>Slab pellet *($F_P=1$)*</u>

$$a_1 = 1 + \frac{M(y-1)}{\coth(My_m) - M(y_m - 1)}$$

$$a_2 = \frac{\cosh(My)}{\cos(My_m) - M(y_m - 1)\sinh(My_m)}$$

$$\theta_V = 2 + M^2(1 - y_m)^2 + 2M(1 - y_m)\tanh(My_m)$$

$$b_2 = \left(1 - \frac{\cosh(My)}{\cosh(My_m)}\right)^2$$

        <u>Spherical pellet ($F_P=3$)</u>

$$a_1 = \left(\frac{1}{1 + (1 - y_m)[My_m \coth(My_m) - 1]}\right)\frac{y_m(1-y)}{y(1 - y_m)} + \frac{1 - y_m/y}{1 - y_m}$$

$$a_2 = a_m \frac{y_m \sinh(My)}{y \sinh(My_m)}$$

$$\theta_V = 2 + \frac{M^2}{3}(1 - y_m)^2(1 + 2y_m) + 2(1 - y_m)[My_m \coth(My_m) - 1]$$

$$b_2 = \left(1 - \frac{y_m \sinh(My)}{y \sinh(My_m)}\right)^2$$

Table 3. The analytical solution for small increments of time in the quantized technique for half order volume reaction model (unsteady state) [96].

**Half order volume reaction model (n=1/2)**

$M = \phi_V b^{1/2}$

**Unsteady state**

First stage

*Slab pellet ($F_P=1$)*

$$a = \frac{\cosh(My)}{\cosh(M)} + \sum_{k=1}^{\infty} \frac{(k\pi)(-1)^k \cos(k\pi y/2)}{(M^2 + k^2\pi^2/4)} \exp\left[-\frac{(M^2 + (k^2\pi^2/4))\theta_V}{\psi\phi_V^2}\right]$$

$$b = \left[1 - \frac{\cosh(My)}{\cosh(My_m)} \frac{\theta_V}{2} - \sum_{k=1}^{\infty} \frac{\frac{(k\pi/2)(-1)^k \cos(k\pi y/2)}{M^2 + (k^2\pi^2/4)} \frac{\psi\phi_V^2}{M^2 + (k^2\pi^2/4)}}{\times \left(1 - \exp\left[-(M^2 + \frac{k^2\pi^2}{4})\frac{\theta_V}{\psi\phi_V^2}\right]\right)}\right]^2$$

*Spherical pellet ($F_P=3$)*

$$a = \frac{\sinh(My)}{y\sinh(M)} + \sum_{k=1}^{\infty} \frac{2k\pi(-1)^k \sin(k\pi y)}{(M^2 + k^2\pi^2)y} \exp\left[-\frac{(M^2 + k^2\pi^2)\theta_V}{\psi\phi_V^2}\right]$$

$$b = \left[1 - \frac{\sinh(My)}{y\sinh(M)} \frac{\theta_V}{2} - \sum_{k=1}^{\infty} \frac{\frac{k\pi(-1)^k \sin(k\pi y)}{(M^2 + k^2\pi^2)y} \frac{\psi\phi_V^2}{M^2 + k^2\pi^2}}{\times \left(1 - \exp\left[-\frac{(M^2 + k^2\pi^2)\theta_V}{\psi\phi_V^2}\right]\right)}\right]^2$$

Second stage ($\theta > \theta_c = 2$)

*Slab pellet ($F_P=1$)*

$\theta_V = 2 + M^2(1 - y_m)^2 + 2M(1 - y_m)\tanh(My_m)$

$$b = \left[1 - \frac{\cosh(My)}{\cosh(My_m)} - \sum_{k=1}^{\infty} \frac{\frac{(k\pi/2y_m)(-1)^k \cos(k\pi y/2y_m)}{M^2 + (k^2\pi^2/4y_m^2)} \frac{\psi\phi_V^2}{M^2 + (k^2\pi^2/4y_m^2)}}{\times \left(1 - \exp\left[-(M^2 + \frac{k^2\pi^2}{4y_m^2})\frac{\theta_c}{\psi\phi_V^2}\right]\right)}\right]^2$$

*Spherical pellet ($F_P=3$)*

$\theta_V = 2 + \frac{M^2}{3}(1 - y_m)^2(1 + 2y_m) + 2(1 - y_m)[My_m \coth(My_m) - 1]$

$$b = \left[1 - \frac{y_m \sinh(My)}{y\sinh(My_m)} - \sum_{k=1}^{\infty} \frac{\frac{(k\pi/y_m)(-1)^k \sin(k\pi y/y_m)}{(M^2 + (k^2\pi^2/y_m^2))y} \frac{\psi\phi_V^2}{M^2 + (k^2\pi^2/y_m^2)}}{\times \left(1 - \exp\left[-(M^2 + \frac{k^2\pi^2}{y_m^2})\frac{\theta_c}{\psi\phi_V^2}\right]\right)}\right]^2$$

Fig. (1) shows a comparison of QM results [92] with approximate solution of Ramachandran [88] and orthogonal collocation solution [40] for steady state half order volume reaction model and the spherical pellet.

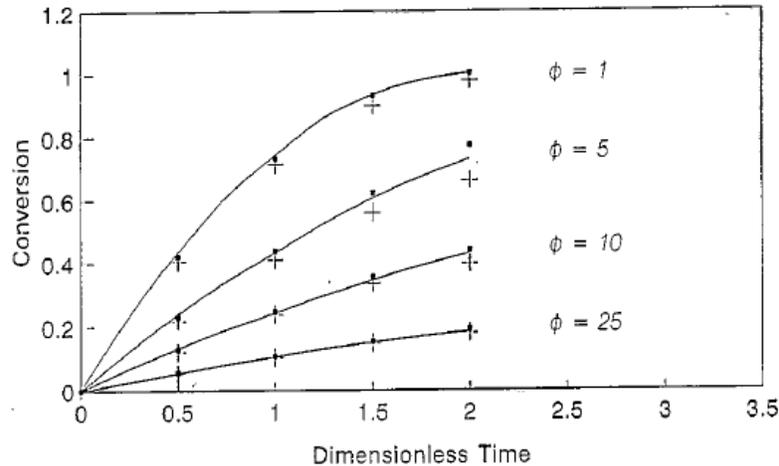

Fig. 1. Comparison of QM results (·) with the approximate solution of Ramachandran (+) and orthogonal collocation solution (solid line) for steady state half order volume reaction model and the spherical pellet [92].

## 2.2. Grain model

Most of the important fluid-solid reactions in chemical and metallurgical industries consist of solid pellets produced from small particles or grains. These pellets are often initially porous, or become porous in the reaction condition. In this model the solid pellet is visualized as consisting of a number of small particles or grains. Surrounding these grains are macropores through which the fluid has to diffuse to reach various grains. The reaction occurs at the surface of each grain according to the sharp interface model. A product layer will form with time in the outer regions of each grain and these will, in turn, offer some resistance to diffusion [46-62].

### 2.2.1. Simple grain model

In this model, the solid product around each grain is highly porous. Thus, the resistance of this product layer is negligible. Considering this point and general assumptions, the dimensionless governing equations in general form are as follows [58]:

$$\psi\sigma^2 \frac{\partial a}{\partial \theta_g} = \frac{\partial^2 a}{\partial y^2} + \frac{F_P - 1}{y}\frac{\partial a}{\partial y} - \sigma^2 r^{*F_g - 1} a \tag{9}$$

$$\frac{\partial r^*}{\partial \theta_g} = -a \tag{10}$$

With the initial and boundary conditions:

$$y = 0; \quad \frac{\partial a}{\partial y} = 0 \tag{11}$$

$$y = 1; \quad a = sh(1 - a) \tag{12}$$

$$\theta_g = 0; \quad r^* = 1 \tag{13}$$

Also, conversion of solid at each time can be calculated by the following equation:

$$X(\theta) = 1 - F_P \int_0^1 y^{F_P - 1} r^{*F_g}(\theta, y) dy \tag{14}$$

This model is studied in two stages. In the first stage ($\theta < \theta_c$), diffusion of gas "A" and reaction between gas "A" and solid "B" are happening simultaneously. At time $\theta = \theta_c$ all solid of the outer layer of the pellet has been reacted. Then second stage is started when $\theta > \theta_c$. In this stage, the gas diffuses through the completely reacted outer layer of the pellet, in order to reach the diffusion-reaction zone. Some approximate [5, 50, 58, 74, 79, 88] and numerical solutions [29, 40, 49, 50, 58] exist in the literature for grain models. The analytical solution for small increments of time in the quantized technique for simple grain model has presented in Table (4) and (5). Also, Fig. (2) shows comparison of QM with numerical and approximate solution of Evan and Ranade [50] for quasi-steady state simple grain model ($F_P=3$, $F_g=2$).

Table 4. The analytical solution for small increments of time in the quantized technique for simple grain model

(quasi-steady state) [91].

**Simple grain model**

$M = \sqrt{\sigma^2 r^{*F_g-1}}$

**Quasi-steady state**

*First stage*

<u>Slab pellet ($F_P=1$) and negligible external mass transfer</u>

$a = \dfrac{\cosh(My)}{\cosh(M)}$, $r^* = 1 - a\theta_g$, $\theta_c = 1$

<u>Spherical pellet ($F_P=3$)</u>

$a = \dfrac{1}{\theta_c} \dfrac{\sinh(My)}{y\sinh(M)}$, $\theta_c = 1 + \dfrac{1}{sh}[M\coth(M) - 1]$, $r^* = 1 - a\theta_g$

*Second stage ($\theta > \theta_c$)*

<u>Slab pellet ($F_P=1$)</u>

$a_1 = 1 + \dfrac{M(y-1)}{\coth(My_m) - M(y_m - 1)}$, $a_2 = \dfrac{\cosh(My)}{\cosh h(My_m) - M(y_m - 1)\sinh(My_m)}$

$\theta_g = 1 + \dfrac{M^2}{2}(1 - y_m)^2 + M(1 - y_m)\tanh(My_m)$

$r_2^* = 1 - \dfrac{\cosh(My)}{\cosh(My_m)}$

<u>Spherical pellet ($F_P=3$)</u>

$a_1 = \left(\dfrac{1}{1 + (1 - y_m + y_m/sh)[My_m\coth(My_m) - 1]}\right)\left(\dfrac{1 - y + y/sh}{1 - y_m + y_m/sh}\right)\dfrac{y_m}{y}$

$+ \dfrac{1 - y_m/y}{1 - y_m + y_m/sh}$

$a_2 = a_m \dfrac{y_m \sinh(My)}{y \sinh(My_m)}$

$\theta_g = 1 + \dfrac{M^2}{6}(1 - y_m)^2(1 + 2y_m) + (1 - y_m)[My_m\coth(My_m) - 1]$

$+ \dfrac{1}{sh}\dfrac{M^2}{3}(1 - y_m^3) + \dfrac{y_m}{N_{sh}}[My_m\coth(My_m) - 1]$

$r_2^* = 1 - \dfrac{y_m \sinh(My)}{y \sinh(My_m)}$

Table 5. The analytical solution for small increments of time in the quantized technique for simple grain model (unsteady state and negligible external mass transfer resistance) [96].

**Simple grain model**

$M = \sqrt{\sigma^2 r^{*F_g - 1}}$

**Unsteady state**

*First stage*

<u>Slab pellet ($F_P=1$)</u>

$$a = \frac{\cosh(My)}{\cosh(M)} + \sum_{k=1}^{\infty} \frac{k\pi(-1)^k \cos(k\pi y/2)}{M^2 + (k^2\pi^2/4)} \exp\left[-\frac{(M^2 + (k^2\pi^2/4))\theta_g}{\psi\sigma^2}\right]$$

$$r^* = 1 - \frac{\cosh(My)}{\cosh(M)}\theta_g$$

$$-\sum_{k=1}^{\infty} \frac{k\pi(-1)^k \cos(k\pi y/2)}{M^2 + (k^2\pi^2/4)}\left(\frac{\psi\sigma^2}{M^2 + (k^2\pi^2/4)}\right)\left(1 - \exp\left[-\frac{(M^2 + (k^2\pi^2/4))\theta_g}{\psi\sigma^2}\right]\right)$$

<u>Spherical pellet ($F_P=3$)</u>

$$a = \frac{\sinh(My)}{y\sinh(M)} + \sum_{k=1}^{\infty} \frac{2k\pi(-1)^k \sin(k\pi y)}{(M^2 + k^2\pi^2)y} \exp\left[-\frac{(M^2 + k^2\pi^2)\theta_g}{\psi\sigma^2}\right]$$

$$r^* = 1 - \frac{\sinh(My)}{y\sinh(M)}\theta_g - \sum_{k=1}^{\infty} \frac{\frac{2k\pi(-1)^k \sin(k\pi y)}{(M^2 + k^2\pi^2)y} \frac{\psi\sigma^2}{M^2 + k^2\pi^2}}{\times\left(1 - \exp\left[-\frac{(M^2 + k^2\pi^2)\theta_g}{\psi\sigma^2}\right]\right)}$$

*Second stage ($\theta > \theta_c$)*

<u>Slab pellet ($F_P=1$)</u>

$$\theta_g = 1 + \frac{M^2}{2}(1 - y_m)^2 + M(1 - y_m)\tanh(My_m)$$

$$r^* = 1 - \frac{\cosh(My)}{\cosh(My_m)}$$

$$-\sum_{k=1}^{\infty} \frac{(k\pi/y_m)(-1)^k \cos(k\pi y/2y_m)}{M^2 + (k^2\pi^2/4y_m^2)}\left(\frac{\psi\sigma^2}{M^2 + (k^2\pi^2/4y_m^2)}\right)\left(1 - \exp\left[-(M^2 + \frac{k^2\pi^2}{4y_m^2})\frac{\theta_c}{\psi\sigma^2}\right]\right)$$

<u>Spherical pellet ($F_P=3$)</u>

$$\theta_g = 1 + \frac{M^2}{6}(1 - y_m)^2(1 + 2y_m) + (1 - y_m)[My_m\coth(My_m) - 1]$$

$$r^* = 1 - \frac{\sinh(My)}{y\sinh(My_m)}\theta_g$$

$$-\sum_{k=1}^{\infty} \frac{(2k\pi/y_m)(-1)^k \sin(k\pi y/y_m)}{(M^2 + (k^2\pi^2/y_m^2))y}\frac{\psi\sigma^2}{M^2 + (k^2\pi^2/y_m^2)}\left(1 - \exp\left[-\left(M^2 + \frac{k^2\pi^2}{y_m^2}\right)\frac{\theta_c}{\psi\sigma^2}\right]\right)$$

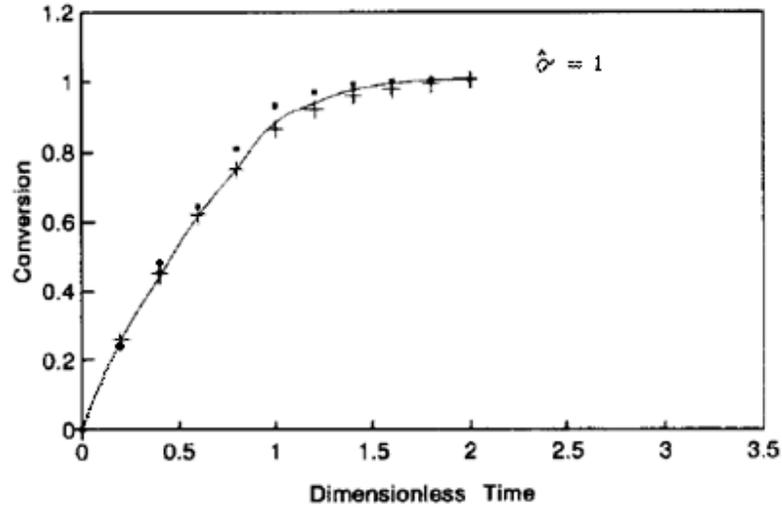

Fig. 2. Comparison of QM (·) with numerical (solid line) and approximate solution (+) of Evans and Ranade for quasi-steady state simple grain model [91].

### 2.2.2. Grain model with product layer resistance

By considering the grain model with product layer resistance, the following equations obtained [59]:

$$\psi\sigma^2 \frac{\partial a}{\partial \theta_g} = \frac{\partial^2 a}{\partial y^2} + \frac{F_P - 1}{y}\frac{\partial a}{\partial y} - \frac{\sigma^2 r^{*F_g} a}{1 - \sigma_g^2 \cdot q'_{F_g}(r^*)} \quad (15)$$

$$\frac{\partial r^*}{\partial \theta_g} = -\frac{a}{1 - \sigma_g^2 \cdot q'_{F_g}(r^*)} \quad (16)$$

With Boundary and initial conditions:

$$y = 0; \quad \frac{\partial a}{\partial y} = 0 \quad (17)$$

$$y = 1; \quad a = 1 \quad (18)$$

$$\theta_g = 0; \quad r^* = 1 \quad (19)$$

And it is known that:

$$q'_{F_g}(r^*) = \frac{\partial q_{F_g}(r^*)}{\partial r^*} \quad (20)$$

For the spherical pellet:

$$q'_{F_g}(r^*) = 1 - 3r^{*2} + 2r^{*3} \tag{21}$$

Where $q_{F_g}(r^*)$ describes progress of the reaction of a grain, under diffusional effects.

The analytical solution for small increments of time in the quantized technique for grain model with product layer resistance is presented in Table (6).

Table 6. The analytical solution for small increments of time in the quantized technique for the grain model with product layer resistance (unsteady state) [96].

**Grain model with product layer resistance**

$$M = \sigma \sqrt{\frac{r^{*2}}{1 + 6\sigma_g^2(r^* - r^{*2})}}$$

**Unsteady state**

*Spherical pellet ($F_P=3$) and grain($F_g=3$)*

First stage

$$a = \frac{\sinh(My)}{y\sinh(M)} + \sum_{k=1}^{\infty} \frac{2k\pi(-1)^k \sin(k\pi y)}{(M^2 + k^2\pi^2)y} \exp\left[-\frac{(M^2 + k^2\pi^2)\theta_g}{\psi\sigma^2}\right]$$

$$r^* + 3\sigma_g^2 r^{*2} - 2\sigma_g^2 r^{*3} = 1 + \sigma_g^2 - \frac{y_m \sinh(My)}{y\sinh(My_m)}\theta_g - \sum_{k=1}^{\infty} \frac{(2k\pi)(-1)^k \sin(k\pi y)}{(M^2 + k^2\pi^2)y} \cdot \frac{\psi\sigma^2}{M^2 + k^2\pi^2} \times \left(1 - \exp\left[-\left(M^2 + \frac{k^2\pi^2}{y_m^2}\right)\frac{\theta_c}{\psi\sigma^2}\right]\right)$$

Second stage ($\theta_g > \theta_c = 1 + \sigma_g^2$)

$$\frac{\theta_g}{1 + \sigma_g^2} = 1 + \frac{M^2}{6}(1 - y_m)^2(1 + 2y_m) + (1 - y_m)[My_m \coth(My_m) - 1]$$

$$r^* + 3\sigma_g^2 r^{*2} - 2\sigma_g^2 r^{*3} = 1 + \sigma_g^2 - \frac{y_m \sinh(My)}{y\sinh(My_m)}\theta_c$$

$$-\sum_{k=1}^{\infty} \frac{(2k\pi / y_m)(-1)^k \sin(k\pi y / y_m)}{(M^2 + (k^2\pi^2 / y_m^2))y} \cdot \frac{\psi\sigma^2}{M^2 + (k^2\pi^2 / y_m^2)} \times \left(1 - \exp\left[-\left(M^2 + \frac{k^2\pi^2}{y_m^2}\right)\frac{\theta_c}{\psi\sigma^2}\right]\right)$$

### 2.2.3. Modified grain model

In most of fluid-solid reactions, significant structural changes will occur due to the difference between molar volumes of the solid product and solid reactant. Thus during the reaction, porosity of the pellet changes and consequently diffusional resistances change.

The modified grain model has been developed for modeling of such reactions [47, 48, 49, 55]. In this model product layer diffusion around each grain as well as variable intergrain diffusion has been considered.

The dimensionless mass balance equations for spherical pellet ($F_P$=3) and grain ($F_g$=3) are as follows [49]:

$$\psi\sigma^2 \frac{\partial a}{\partial \theta_g} = \frac{1}{y^2}\frac{\partial}{\partial y}\left(\delta y^2 \frac{\partial a}{\partial y}\right) - \frac{\sigma^2 r^{*2} a}{1+6\sigma_g^2(r^* - \frac{r^{*2}}{r^{**}}))} \tag{22}$$

$$\frac{\partial r^*}{\partial \theta_g} = -\frac{a}{1+6\sigma_g^2(r^* - \frac{r^{*2}}{r^{**}})} \tag{23}$$

With boundary and initial conditions:

$$y = 0; \quad \frac{\partial a}{\partial y} = 0 \tag{24}$$

$$y = 1; \quad a = 1 \tag{25}$$

$$\theta_g = 0; \quad r^* = 1 \tag{26}$$

For spherical grains the dimensionless outer radius of grain is:

$$r^{**} = [Z_\upsilon + (1-Z_\upsilon)r^{*3}]^{1/3} \tag{27}$$

Where

$$Z_\upsilon = \frac{\nu_D \rho_B M_D}{\nu_B \rho_D M_B (1-\varepsilon_D)} \tag{28}$$

The porosity of the pellet at each time can be related to the grain size:

$$\frac{1-\varepsilon}{1-\varepsilon_0} = r^{**3} \tag{29}$$

Also, the ratio of the intergrain diffusion at each time to time zero is as follows:

$$\delta(r^*) = \frac{D_e}{D_{eo}} = \left(\frac{\varepsilon}{\varepsilon_0}\right)^2 = \left[1 - \frac{1-\varepsilon_0}{\varepsilon_0}(Z_\upsilon - 1)(1 - r^{*3})\right]^2 \tag{30}$$

The analytical solution for small increments of time in the quantized technique for modified grain model is presented in Tables (7) and (8). Also, Fig. (3) shows the comparison of QM with the orthogonal collocation solution of Dudukovic and Lamba [40], and the approximate solution of Ramachandran [88] for the first-order volume reaction model and a spherical pellet. As Fig. (3) shows, there is good agreement between QM prediction and existing numerical solutions.

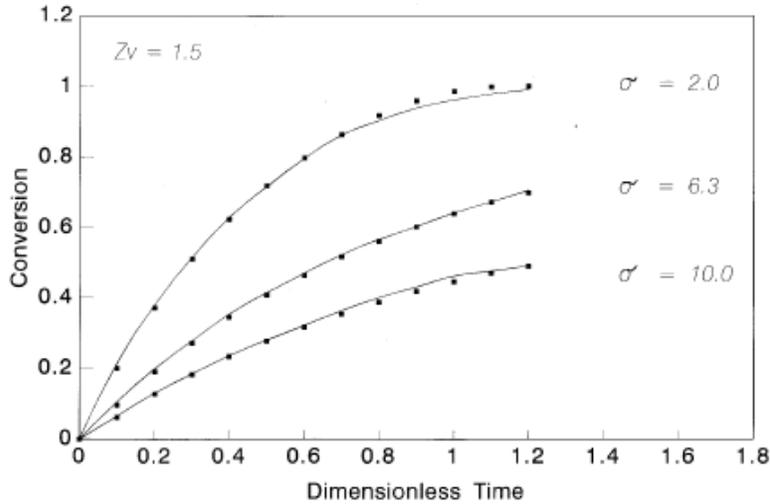

Fig. 3. Comparison of QM (•) with the orthogonal collocation solution (solid line) of Prasannan et al. , for the modified grain model with Zv=1.5, $\sigma_g^2$=0.167, $\varepsilon_0$=0.5, and some values of σ [94].

Table 7. The analytical solution for small increments of time in the quantized technique for the modified grain model (quasi-steady state) [94].

**Modified grain model**

$$M = \sqrt{\frac{\sigma^2 r^{*2}}{\left[1+6\sigma_g^2\left(r^* - \frac{r^{*2}}{[Z_v + (1-Z_v)r^{*3}]^{1/3}}\right)\right].\delta(r^*)}}$$

**Quasi-steady state**

*Spherical pellet ($F_P=3$) & Spherical grain ($F_g=3$)*

First stage

$$a = \frac{\sinh(My)}{y\sinh(M)}$$

$$r^* + 3\sigma_g^2 r^{*2} + \frac{3\sigma_g^2}{Z_v - 1}[Z_v + (1-Z_v)r^{*3}]^{2/3} = 1 + (3 + \frac{3}{Z_v - 1})\sigma_g^2 - a\theta_g$$

Second stage ($\theta_g > \theta_c = 1 + 3\sigma_g^2 + \frac{3\sigma_g^2}{Z_v - 1}(1 - Z_v^{2/3})$

$$a_1 = \left(\frac{1}{1 + (1-y_m)[My_m \coth(My_m) - 1]}\right)\left(\frac{1-y}{1-y_m}\right)\frac{y_m}{y} + \frac{1 - y_m/y}{1 - y_m}$$

$$a_2 = \left(\frac{1}{1 + (1-y_m)[My_m \coth(My_m) - 1]}\right)\frac{y_m \sinh(My)}{y \sinh(My_m)}$$

$$r^* + 3\sigma_g^2 r^{*2} + \frac{3\sigma_g^2}{Z_v - 1}[Z_v + (1-Z_v)r^{*3}]^{2/3}$$

$$= (1 + 3\sigma_g^2 + \frac{3\sigma_g^2}{Z_v - 1}) - (1 + 3\sigma_g^2 + \frac{3\sigma_g^2}{Z_v - 1}[1 - Z_v^{2/3}])\frac{y_m \sinh(My)}{y \sinh(My_m)}$$

$$\frac{\theta_g}{1 + 3\sigma_g^2 + \frac{3\sigma_g^2}{Z_v - 1}(1 - Z_v^{2/3})}$$

$$= 1 + \frac{M^2}{6}(1 - y_m)(1 + 2y_m) + (1 - y_m)[My_m \coth(My_m) - 1]$$

Table 8. The analytical solution for small increments of time in the quantized technique for the modified grain model (unsteady state) [96].

---

**Modified grain model**

$$M = \sigma \sqrt{\frac{r^{*2}}{[1+6\sigma_g^2(r^* - (r^{*2}/r^{**}))]\delta(r^*)}}$$

**Unsteady state**

*Spherical pellet ($F_P=3$) & Spherical grain ($F_g=3$)*

First stage

$$a = \frac{\sinh(My)}{y\sinh(M)} + \sum_{k=1}^{\infty} \frac{2k\pi(-1)^k \sin(k\pi y)}{(M^2+k^2\pi^2)y} \exp\left[-\frac{(M^2+k^2\pi^2)\theta_g}{\psi\sigma^2}\right]$$

$$r^* + 3\sigma_g^2 r^{*2} - \frac{3\sigma_g^2}{Z-1}[Z+(1-Z)r^{*3}]^{2/3} =$$

$$1+(3+\frac{3}{Z-1})\sigma_g^2 - \frac{\sinh(My)}{y\sinh(M)}\theta_g - \sum_{k=1}^{\infty} \frac{(2k\pi)(-1)^k \sin(k\pi y)}{(M^2+k^2\pi^2)y} \frac{\psi\sigma^2/\delta}{M^2+k^2\pi^2} \times \left(1-\exp\left[-\frac{(M^2+k^2\pi^2)\theta_g}{\psi\sigma^2/\delta}\right]\right)$$

Second stage ($\theta_g > \theta_c = 1+3\sigma_g^2 + (3\sigma_g^2/(Z-1))(1-Z^{2/3})$)

$$\frac{\theta_g}{\theta_c} = 1 + \frac{M^2}{6}(1-y_m)^2(1+2y_m) + (1-y_m)[My_m \coth(My_m) - 1]$$

$$r^* + 3\sigma_g^2 r^{*2} + \frac{3\sigma_g^2}{Z-1}[Z+(Z-1)r^{*3}]^{2/3} = 1 + \left(3+\frac{3}{Z-1}\right)\sigma_g^2 - \frac{y_m \sinh(My)}{y\sinh(My_m)}\theta_c$$

$$-\sum_{k=1}^{\infty} \frac{(2k\pi/y_m)(-1)^k \sin(k\pi y/y_m)}{(M^2+(k^2\pi^2/y_m^2))y} \frac{\psi\sigma^2/\delta}{M^2+(k^2\pi^2/y_m^2)}\left(1-\exp\left[-\left(M^2+\frac{k^2\pi^2}{y_m^2}\right)\frac{\theta_c}{(\psi\sigma^2/\delta)}\right]\right)$$

---

## 2.3. Random pore model

In most of fluid–solid reactions, the complicated pore size distribution exists in the solid pellet. Moreover, significant structural changes occur due to the difference between molar volumes of the solid product and solid reactant. Therefore, the size of pores and diffusional resistances change by the progress of the reaction.

The random pore model is one of the most sophisticated approaches in the fluid–solid reactions. In this model, product layer diffusion around each pore as well as variable diffusion in the cylindrical pores was considered [63-65]. Random pore model has been used for modeling the CaO + SO$_2$ reaction [70]. This model was also modified to consider of the maximum internal surface area in the gasification reactions [69]. Assumptions of this model are similar to assumption of modified grain model.

In this case, the dimensionless mass conservation equations within spherical particle are as follows [98]:

$$\psi \cdot \phi_r^2 \frac{\partial a}{\partial \theta_r} = \frac{1}{y^2} \frac{\partial}{\partial y}(\delta y^2 \frac{\partial a}{\partial y}) - \frac{\phi_r^2 ab\sqrt{1-\Psi \ln(b)}}{1+(\beta Z/\Psi)[\sqrt{1-\Psi \ln(b)}-1]} \quad (31)$$

$$\frac{\partial b}{\partial \theta_r} = -\frac{a.b.\sqrt{1-\Psi \ln(b)}}{1+\frac{\beta.Z}{\Psi}(\sqrt{1-\Psi \ln(b)}-1)} \quad (32)$$

$$y = 0; \quad \frac{\partial a}{\partial y} = 0 \quad (33)$$

$$y = 1; \quad \frac{\partial a}{\partial y} = (sh/\delta)(1-a) \quad (34)$$

$$\theta_r = 0; \quad b = 1 \quad (35)$$

Where Z is the ratio of molar volumes of the solid product to the solid reactant:

$$Z = \frac{v_D \rho_B M_D}{v_B \rho_D M_B} \quad (36)$$

$$\delta = \left(\frac{\varepsilon}{\varepsilon_0}\right)^2 = \left[1 - \frac{(Z-1)(1-\varepsilon_0)(1-b)}{\varepsilon_0}\right]^2 \quad (37)$$

Also, conversion of solid at each time can be calculated by the following Eq. (7).

Some approximate [75, 87] and numerical solutions [64, 87] exist in the literature for grain models. The analytical solution for small increments of time in quantized technique for random pore model is presented in Table (9). Furthermore, the comparison for quasi-steady state random pore model ($F_P=3$, $F_g=2$) with Ref. [64] is presented in Fig. (4) successfully. As these results show, there is a good agreement between QM predictions and numerical solution especially for low and intermediate values of Thiele modulus.

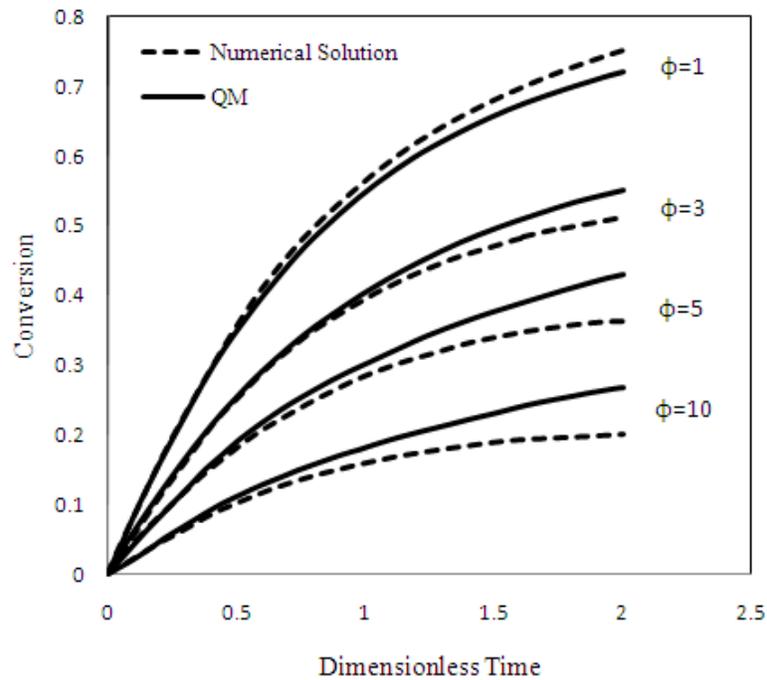

Fig. 4. Comparison of QM with numerical solution of Bhatia and Perlmutter [65] for quasi-steady state random pore model.

Table 9. The analytical solution for small increments of time in the quantized technique for the random pore model.

**Random pore model**

**Quasi-steady state**

$$M = \frac{\phi_r \sqrt{\dfrac{b.\sqrt{1-\Psi \ln(b)}}{1+\dfrac{\beta.Z}{\Psi}(\sqrt{1-\Psi \ln(b)}-1)}}}{\left[1-\dfrac{(Z-1)(1-\varepsilon_0)(1-b)}{\varepsilon_0}\right]}$$

*Spherical pellet ($F_P=3$)*

$$a = \left(\frac{1}{1+\dfrac{\delta}{sh}[M\cot(M)-1]}\right)\left(\frac{\sinh(My)}{y\sinh(M)}\right)$$

$$\frac{2}{\psi}(1-\frac{\beta Z}{\psi})(1-\sqrt{1-\psi \ln b}) + \frac{\beta Z}{\psi}\ln b = \frac{\theta_r}{\theta_c}\frac{\sinh(My)}{y\sinh(M)}$$

**Unsteady state & β=0 [96]**

$$M = \phi_r \sqrt{b\sqrt{1-\Psi \ln(b)}}$$

*Spherical pellet ($F_P=3$)*

$$a = \frac{\sinh(My)}{y\sinh(M)} + \sum_{k=1}^{\infty}\frac{2k\pi(-1)^k \sin(k\pi y)}{(M^2+k^2\pi^2)y}\exp\left[-\frac{(M^2+k^2\pi^2)\theta_g}{\psi\phi_r}\right]$$

$$\Psi \ln(b) = 1 - \left\{1+\frac{\Psi}{2}\left[\frac{\sinh(My)}{y\sinh(M)}\theta_g + \sum_{k=1}^{\infty}\frac{2k\pi(-1)^k \sin(k\pi y)}{(M^2+k^2\pi^2)y}\frac{\psi\phi_r^2}{M^2+k^2\pi^2}\times\left(1-\exp\left[-\frac{(M^2+k^2\pi^2)\theta_g}{\psi\phi_r^2}\right]\right)\right]\right\}^2$$

## 2.4. Nucleation model

Most gas-solid reactions start with formation of nuclei at the solid surface. As the reaction progresses, this nuclei increase in size, overlap with one another and cover the surface. Nucleation effects are often significant, for example in reduction of metallic oxides. A more general problem has been presented for isothermal and non-isothermal cases [68]. The effect of pore diffusion in the nucleation model has been analyzed elsewhere [74].

Assumptions of this model are similar to the assumption of volume reaction model. The general dimensionless conservation equations of fluid and solid, based on nucleation growth kinetics are as follows [74]:

$$\psi \sigma_N^2 \frac{\partial a}{\partial \theta} = \frac{\partial^2 a}{\partial y^2} + \frac{F_p - 1}{y} \frac{\partial a}{\partial y} + 2F_P \sigma_N^2 \frac{a}{f'(b)} \quad (38)$$

$$\frac{\partial b}{\partial \theta_n} = \frac{a}{f'(b)} \quad (39)$$

With the initial and boundary conditions:

$$y = 0; \quad \frac{\partial a}{\partial y} = 0 \quad (40)$$

$$y = 1; \quad a = 1 \quad (41)$$

$$\theta_n = 0; \quad b = 1$$

Where the dimensionless quantities are defined in the notation, and

$$f'(b) = \frac{\partial f(b)}{\partial b} \quad (42)$$

$$f(b) = [-\ln(b)]^{1/n} \quad (43)$$

This problem has been solved for two different values of n (n=1, n=3). By considering n=1, the nucleation model reduces to the first order volume reaction model. Conversion of solid at each time can be calculated by the following Eq. (7).

Some approximate [25, 88] and numerical solutions [81, 40, 99] exist in the literature for grain models. The analytical solution for small increments of time in the quantized technique for the nucleation model (n=3) is presented in Table (10). Also, Figs. (5) and (6) display a comparison of QM results with the numerical and approximate solution of Sohn and Kim [75] for the nucleation model. According to Figs. (5) and (6) and other data that exist in Jamshidi and Ale Ebrahim's [93] work, the agreement between numerical solution and QM is excellent for n=1, and good for n=3, at small modulus. For large modulus and small and intermediate dimensionless times, QM is more accurate than Sohn's approximate solution. However, for large modulus and large dimensionless times, Sohn's approximate solution gives better results [93].

Table 10. The analytical solution for small increments of time in the quantized technique for the nucleation model (n=3) [93, 96].

**Nucleation model**

$$M = \sigma_N \sqrt{2F_P / f'(b)}$$

**Quasi-steady state**

*Spherical pellet ($F_P=3$)*

$$a = \frac{\sinh(My)}{y \sinh(M)}$$

$$b = \exp[(-a\theta_n)]^3$$

**Unsteady state**

*Spherical pellet ($F_P=3$)*

$$a = \frac{\sinh(My)}{y \sinh(M)} + \sum_{k=1}^{\infty} \frac{2k\pi(-1)^k \sin(k\pi y)}{(M^2 + k^2\pi^2)y} \exp[\frac{-(M^2 + k^2\pi^2)\theta_n}{\psi \sigma_N}]$$

$$b = \exp\left[-\left[\frac{\sinh(My)}{y \sinh(M)}\theta_g + \sum_{k=1}^{\infty} \frac{2k\pi(-1)^k \sin(k\pi y)}{(M^2 + k^2\pi^2)y} \frac{\psi \sigma_N^2}{M^2 + k^2\pi^2} \times \left(1-\exp\left[-\frac{(M^2 + k^2\pi^2)\theta_g}{\psi \sigma_N^2}\right]\right)\right]^3\right]$$

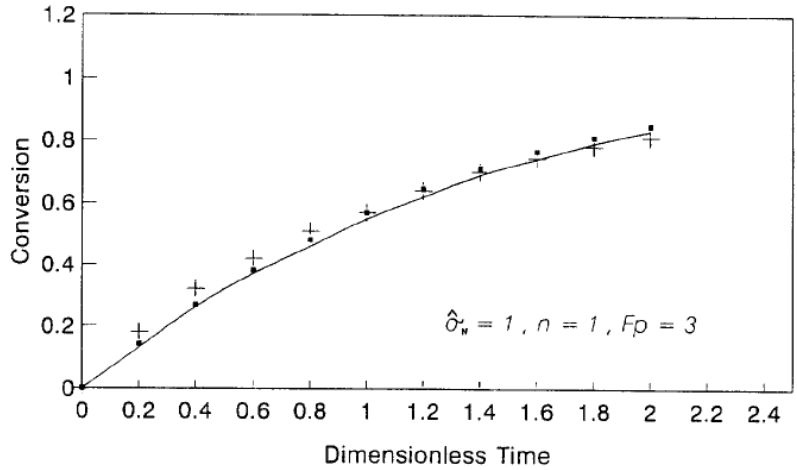

Fig. 5. Comparison of QM results (•) with numerical (solid line) and approximate solution (+) of Sohn for nucleation model, n=1 [93].

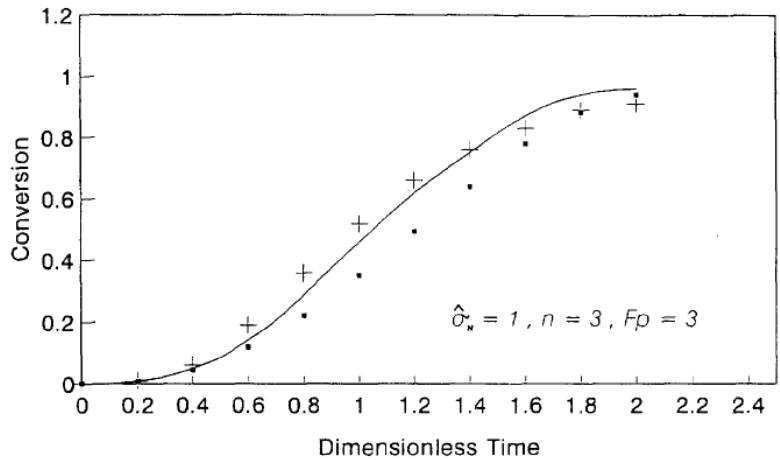

Fig. 6. Comparison of QM results (•) with numerical (solid line) and approximate solution (+) of Sohn for nucleation model, n=3 [93].

**2.5. Simultaneous fluid-solid reaction model**

In many industrial operations, more than one reaction can take place simultaneously. Because the gas (or solid) contains multiple reactant, or the products of reaction are reactive with the solid. Consecutive reaction [100, 101], reaction of two different gases with one solid [85, 102], reaction of two solid with a gas [16, 28, 86, 103], and reaction of gaseous product with the second solid [104, 105] are some commonly cases of simultaneous gas solid reactions.

Wenand and Wei [106], and Rehmat and Saxena [22] studied the problem of simultaneous gas solid reactions on the nonisothermal shrinking core model. However, for pellets of considerable porosity Sohn and Braun [85], Sohn and Rajamani [84], and Tone and Wen [107] presented a more general model.

The reaction between a porous solid and two gaseous reactants can be represented as:

$$A(g) + v_{B1}B(s) \to v_P P_1(g) + v_{D1} D(s) \tag{44}$$

$$C(g) + v_{B2}B(s) \to n_Q P_2(g) + v_{D2} D(s) \tag{45}$$

These reactions are important in the direct reduction of iron oxide by reformed natural gas. The assumptions of this model are similar to the assumptions of volume reaction model.

The general quasi-steady state dimensionless conservation equations of this model are as follows [72]:

$$\frac{\partial^2 \psi_A}{\partial y^2} + \frac{F_p - 1}{y} \frac{\partial \psi_A}{\partial y} = -2 F_P \sigma_A^2 \frac{\psi_A}{f'(b)} \tag{46}$$

$$\frac{\partial^2 \psi_C}{\partial y^2} + \frac{F_p - 1}{y} \frac{\partial \psi_C}{\partial y} = -2 F_P \sigma_C^2 \frac{\psi_C}{f'(b)} \tag{47}$$

With boundary conditions:

$$y = 0; \quad \frac{\partial \psi_A}{\partial y} = \frac{\partial \psi_C}{\partial y} = 0 \tag{48}$$

$$y = 1; \quad \psi_A = \psi_{Ab} \tag{49}$$

$$y = 1; \quad \psi_C = \psi_{Cb} \tag{50}$$

And

$$f'(b) = \frac{\partial f(b)}{\partial b} \tag{51}$$

$$f(b) = [-\ln(b)]^{1/n} \tag{52}$$

For solid reactant describing equation are:

$$\frac{\partial b}{\partial \theta_{tg}} = \frac{\psi_A + \psi_C}{f'(b)} \tag{53}$$

$$\frac{\partial b_A}{\partial \theta_{tg}} = \frac{\psi_A}{f'(b)} \tag{54}$$

With an initial condition:

$$\theta_{tg} = 0; \quad b = b_A = 1 \tag{55}$$

The analytical solution for small increments of time in the quantized technique for the model of reaction between a porous solid and two gaseous reactants is presented in Table (11). The conversion of spherical pellet at each time can be calculated from:

$$X = 1 - 3\int_0^1 b.y^2.dy \tag{56}$$

The conversion due to gas "A" can be calculated from:

$$X_A = 1 - 3\int_0^1 b_A.y^2.dy \tag{57}$$

Finally, overall selectivity of system can be calculated from:

$$S_0 = \frac{X_A}{X - X_A} \tag{58}$$

This solution technique is developed for the case where both the diffusion and chemical reactions are important. The predicted conversion time behavior is compared with the numerical solution of Ref. [85]. Fig. (7) Presents this comparison for a given $\sigma_A$, $\sigma_C$ and $\psi_{Ab}$. Also, Fig. (8) shows comparison of this work with numerical solution of Ref. [85] for overall selectivity. Like other cases, the use of QM presents rapid calculations with a reasonable accuracy.

Table 11. The analytical solution for small increments of time in the quantized technique for model of reaction between a porous solid and two gaseous reactants.

| **Model of reaction between a porous solid and two gaseous reactants (n=1)** |
|---|
| $M_A = \sigma_A \sqrt{2F_P / f'(b)}$ <br> $M_C = \sigma_C \sqrt{2F_P / f'(b)}$ <br><br> *Spherical pellet ($F_P=3$)* <br> $\psi_A = \psi_{Ab} \dfrac{\sinh(M_A y)}{y \sinh(M_A)}$ <br> $\psi_C = \psi_{Cb} \dfrac{\sinh(M_C y)}{y \sinh(M_C)}$ <br> $b = \exp[(-(\psi_A + \psi_B).\theta_{tg})]$ <br> $1 - b_A = \dfrac{\psi_A}{\psi_A + \psi_C}(1 - \exp[-(\psi_A + \psi_C).\theta_{tg}])$ |

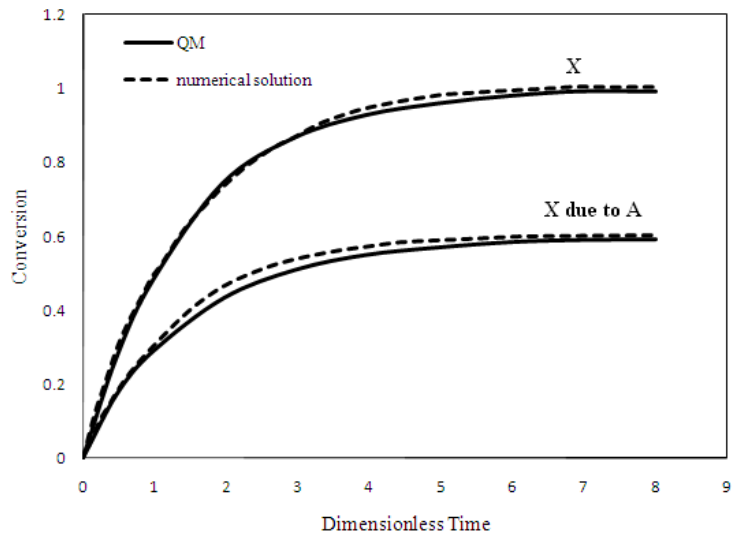

Fig. 7. Comparison of this work (QM) with numerical solution of Ref. [85], for conversion-time behavior of simultaneous reactions at a small value of $\sigma_A$ and a large value of $\sigma_C$ ($\sigma_A=0.1$, $\sigma_C=3$, $\psi_{Ab}=0.4$).

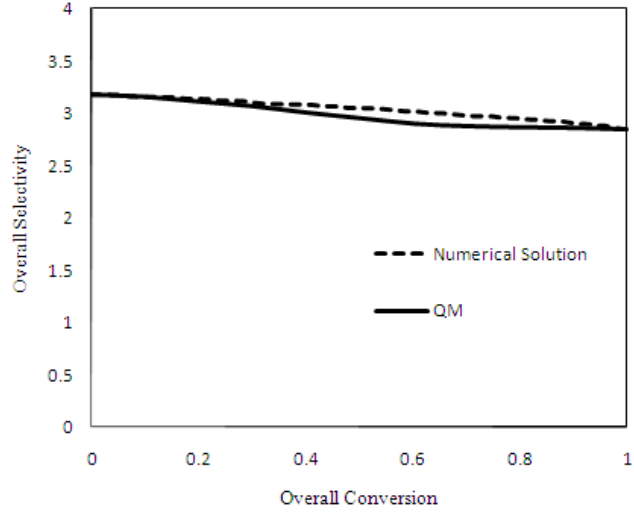

Fig. 8. Comparison of this work (QM) with numerical solution of Ref. [85], for overall selectivity of simultaneous reactions at a small value of $\sigma_A$ and a large value of $\sigma_C$ ($\sigma_A$=0.1, $\sigma_C$=3, $\psi_{Ab}$=0.4).

## 3. Solving PDEs in gas-solid reactions in packed beds

Consider the following chemical reaction has been done between fluid and solid in a packed bed:

$$v_A A(f) + v_B B(s) \rightarrow v_P P(f) + v_G G(s) \tag{59}$$

Where f and s show fluid and solid phases, respectively. Also, rate of reaction, R, in the unit volume of solid is defined as:

$$R = k f(X) C_A \tag{60}$$

Where k is rate constant and the function f(X) can be selected from the different reaction models or with respect to empirical equations.

The dimensionless mass balance equations in addition to their boundary conditions for A and B components have been written as follows [108]:

$$\frac{\partial^2 Y}{\partial \eta^2} - Pe \frac{\partial Y}{\partial \eta} = \beta(Y - a|_{\rho=1}) \tag{61}$$

$$\eta = 0; \quad g(\tau) = Y - \frac{1}{Pe}\frac{\partial Y}{\partial \eta} \tag{62}$$

$$\eta = \Lambda; \quad \frac{\partial Y}{\partial \eta} = 0 \tag{63}$$

$$\frac{1}{y^2}\frac{\partial}{\partial y}\left(y^2\frac{\partial a}{\partial y}\right) = \Phi^2 f(X)a \tag{64}$$

$$y = 0; \quad \frac{\partial a}{\partial y} = 0 \tag{65}$$

$$y = 1; \quad \frac{\partial a}{\partial y} = Bi_m(Y - a) \tag{66}$$

$$\frac{\partial X}{\partial \tau} = f(X)a \tag{67}$$

$$\tau = 0; \quad X = 0 \tag{68}$$

Where dimensionless variables are defined as follows:

$$\eta = \frac{x}{R}, \quad \rho = \frac{r}{R}, \quad a = \frac{C_A}{C_{A0}}, \quad Y = \frac{C_{Ab}}{C_{A0}}, \quad \tau = \nu_B \frac{kC_{A0}t}{\sigma_0} \tag{69}$$

Also, dimensionless parameters are as follows:

$$Pe = \frac{uR}{D_L}, \quad \beta_b = \frac{(1-\varepsilon_b)}{\varepsilon_b}\frac{3k_mR}{D_L}, \quad \Phi^2 = \frac{\nu_A kR}{D_e}, \quad Bi_m = \frac{k_mR}{D_e} \tag{70}$$

The chemical reaction has been considered irreversible and first order with respect to both gas and solid reactants. Therefore Function f(X) can be written as:

$$f(X) = 1 - X \tag{71}$$

Also function g(τ)=1, that means no variation of entrance concentration with time.

Fernandes and Georg [108] proposed a new method for solution of this system. They extended the semianalytical technique of Delborghi et al. [83] and Dudukovic and Lamba [40] for fluid-solid reaction in single pellet to reactions in packed beds. Also, this system has been solved by finite difference method [98].

In quantized method, it was assumed that for small increments of time, X can be treated as a constant in Eq. (64). Therefore, X (j-1, i, k) were used instead of X (j, i, k) in Eq. (64) as an approximation where j, i and k are the time index, the position index for the pellet and the position index for the bed, respectively. Using this technique, a new modified Thiele Modulus is defined for packed beds as:

$$M(j,i,k) = \Phi\sqrt{1 - X(j-1,i,k)} \qquad (72)$$

In this method, the procedure of solution is as follows:

1. Calculation of concentration in surface of the pellet and bottom of the bed, at time zero with attention to conversion of zero in this time based on analytical solution of Eq. (64)

2. Calculation of concentration in bulk of gas based on analytical solution of Eq. (61).

3. Iteration of calculation in axial direction of the bed from the bottom to the top.

4. Calculation of conversion with respect to results of concentration in surface of pellet for total height of the bed.

5. One step time change and start of computation from step 1.

6. Calculation of concentrations in pellets in each section of the bed, with attention to concentration Y given from last 5 steps for total height of the bed for all times.

The analytical solution of Eqs. (61), (64) and (67) for small increments of time in quantized technique is presented in Table (12). Fig. (9) shows comparison of this work (QM) with semi-analytical technique of Ref. [108], for dimensionless gas bulk concentration (Y) versus dimensionless height of the bed at various times. Also, Fig. (10) shows comparison of this work (QM) with semi-analytical technique of Ref. [108] for cumulative concentration versus dimensionless height of the bed at various times. Cumulative concentration is defined as follows:

$$C_Y = \int_0^\tau Y(\eta, \tau') d\tau' \qquad (73)$$

As these figures show, there is a good agreement between QM predictions and solution of Ref. [108] for a packed bed reactor.

Table 12. The analytical solution for small increments of time in the quantized technique for the gas-solid reactions model in packed beds.

**Model of gas-solid reactions in packed beds**

$$M(j,i,k) = \Phi\sqrt{1 - X(j-1,i,k)}$$

$$a = \frac{Bi_m \sinh(My)(Y/y)}{M\cosh(M) + (Bi_m - 1)\sinh(M)}$$

$$Y = \frac{Pe(1 - a|_{y=1})}{Pe - Pe\frac{r_2}{r_1}e^{(r_2 - r_1)\Lambda} - r_2 + r_2^{(r_2 - r_1)\Lambda}}\left[e^{r_2\eta} - \frac{r_2}{r_1}e^{(r_2 - r_1)\Lambda}e^{r_1\eta}\right] + a|_{y=1}$$

$$X = 1 - e^{-a\tau}$$

$$r_{1,2} = \frac{Pe \pm \sqrt{Pe^2 + 4\beta}}{2}$$

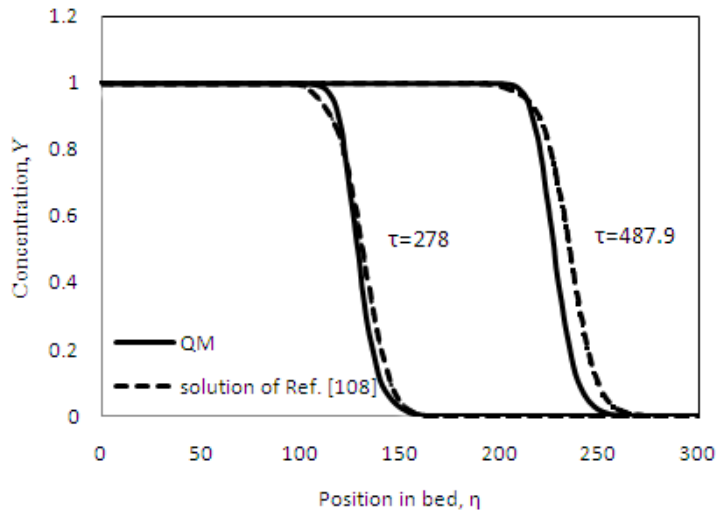

Fig. 9. Comparison of this work (QM) with semi-analytical technique of Ref. [108], for dimensionless bulk gas concentration (Y) versus dimensionless height of the bed at various times, $\Phi=10$, $Bi_m=50$, $Pe=1.1$, $\beta=3.3$.

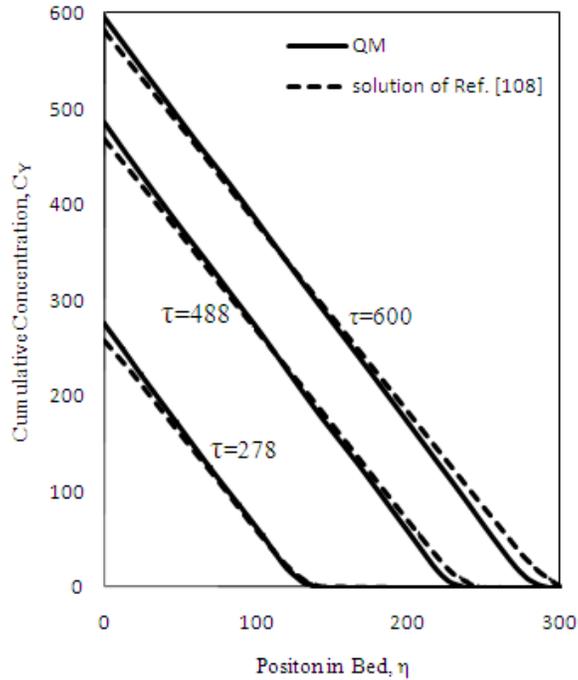

Fig. 10. Comparison of this work (QM) with semi-analytical technique of Ref. [108], for cumulative concentration versus dimensionless height of the bed at various times, $\Phi=10$, $Bi_m=50$, $Pe=1.1$, $\beta=3.3$.

## 4. Conclusion

In this work, the quantized solution was reviewed for various mathematical models in gas-solid and liquid-solid reactions. Furthermore, this technique was applied for quasi-steady state random pore model, simultaneous reactions, and gas-solid reactions in packed beds successfully. The numerical solution of these equations is computationally tedious, and has some problems in the convergence of the results. As this work shows, for low and medium Thiele modulus in a single pellet reaction, there is a good agreement between QM results and other numerical solutions. Also, quantized solution can be applied for the rapid conversion-time predictions and estimation of model parameters from experimental data with the reduced computational effort.

**Nomenclature**

| | |
|---|---|
| $a = C_A/C_{Ab}$ | dimensionless fluid concentration |
| $b = C_B/C_{B0}$ | dimensionless solid concentration |
| $C_A$ | fluid concentration in pellet |
| $C_{Ab}$ | bulk fluid concentration |
| $C_{A0}$ | Initial fluid concentration |
| $C_B$ | solid concentration |
| $C_{B0}$ | initial solid concentration |
| $D_e$ | effective diffusivity in the porous pellet |
| $D_{e0}$ | initial effective diffusivity |
| $D_L$ | axial dispersion coefficient within column |
| $D_p$ | diffusion coefficient in product layer |
| $F_g$ | shape factor of the grain |
| $F_P$ | shape factor of the pellet |
| I | position index for the pellet |
| J | time index |
| K | position index for the bed |
| $k_m$ | mass transfer coefficient |
| $k_s$ | surface rate constant |
| $k_v$ | volume rate constant |
| $k_1, k_2$ | Reaction rate constants for Eqs. (44) and (45) |
| $K_1, K_2$ | Equilibrium constants for Eqs. (44) and (45) |
| M | modified Thiele modulus |
| $M_B$ | molecular weight of solid reactant |
| $M_D$ | molecular rate of solid product |
| N | reaction order with respect to solid |
| R | position of each point in the pellet |
| $r_g$ | outlet radius of the grains |
| $r_{g0}$ | Initial grain radius |
| $r_{gc}$ | radius of unreacted in the grains |
| $r^* = r_{gc}/r_{g0}$ | dimensionless unreacted core radius in the grains |
| $r^* = r_g/r_{g0}$ | dimensionless outlet grain radius |
| R | pellet characterization length |
| $R_b$ | bed diameter |
| $S_0$ | reaction surface per unit volume |
| T | Time |
| X | position of each point in the bed |
| X | solid conversion |
| $y = r/R$ | dimensionless position in pellet |
| $Y = C_{Ab}/C_{A0}$ | dimensionless bulk concentration |
| Z | volume change parameter |

*Greek letter*

| | |
|---|---|
| $\beta = 2k_s(1-\varepsilon_0)/v_B D_P S_0$ | product layer resistance in the random pore model |
| $\Delta$ | dimensionless effective gas diffusivity |
| E | pellet porosity |
| $\varepsilon_0$ | initial pellet porosity |
| $\varepsilon_D$ | porosity of the product layer |
| $\eta = x/R$ | Dimensionless bed length |
| $\theta_g = \nu_B k_s C_{Ab} M_B t / (\rho_B r_{g0})$ | dimensionless time in the grain model |
| $\theta_n = \nu_B k_v C_{Ab} t$ | dimensionless time in the nucleation model |
| $\theta_r = k_s S_0 C_{Ab} t / (C_{B0}(1-\varepsilon_0))$ | dimensionless time in the random pore model |
| $\theta_{tg} = \left[ \nu_{B1} k_1 (C_{Ag} - C_{Pg}/K_1) + \nu_{B2} k_2 (C_{Cg} - C_{Qg}/K_2) \right]$ | Dimensionless time for reaction of two fluid with one solid |
| $\theta_v = \nu_B k_v C_{Ab} C_{B0}^{n-1} t$ | dimensionless time in the volume reaction model |
| $\nu_B$ | stoichiometric coefficient of solid reactant |
| $\nu_D$ | stoichiometric coefficient of solid product |
| $\rho_B$ | density of the solid reactant |
| $\rho_B$ | density of the solid product |
| $\sigma = R\sqrt{F_g k_s (1-\varepsilon_0)/(D_{e0} r_{g0})}$ | Thiele modulus of the grain models |
| $\sigma_g = \sqrt{k_s r_{g0}/(2 F_g D_P)}$ | Thiele modulus of the grains |
| $\sigma_N = R\sqrt{k_v \rho_B (1-\varepsilon)/(2 F_P D_e M_B)}$ | Thiele modulus of the nucleation model |
| $\tau$ | dimensionless time |
| $\psi = \varepsilon C_{Ab}/((1-\varepsilon)C_{B0})$ | accumulation parameter |
| $\Psi$ | random pore model parameter |
| $\psi_A = \nu_{B1} k_1 (C_A - C_P/K_1)/(\nu_{B1} k_1 (C_{Ag} - C_{Pg}/K_2) + \nu_{B2} k_2 (C_{Cg} - C_{Qg}/K_2))$ | dimensionless concentration of fluid "A" |
| $\psi_{Ag}$ | dimensionless concentration of fluid "A" in bulk |
| $\psi_{AC} = \nu_{B2} k_2 (C_C - C_Q/K_2)/(\nu_{B1} k_1 (C_{Ag} - C_{Pg}/K_2) + \nu_{B2} k_2 (C_{Cg} - C_{Qg}/K_2))$ | dimensionless concentration of fluid "C" |
| $\psi_{Cg} = 1 - \psi_{Ag}$ | dimensionless concentration of fluid "C" in bulk |

$$\phi_r = R\sqrt{k_s S_0 /(\nu_B D_{e0})}$$  Thiele modulus of the random pore model

$$\phi_v = R\sqrt{k_v C_{B0}^n / D_{e0})}$$  Thiele modulus of the volume reaction model